\documentclass[leqno,a4paper]{article}

\usepackage[english]{babel}
\usepackage{amsmath}
\usepackage{amsthm}
\usepackage{amssymb}
\usepackage{enumerate}
\usepackage{xspace}
\usepackage{euscript}
\usepackage{graphicx}
\usepackage{amscd}
\usepackage{epsfig}
\usepackage{tabularx}
        \newtheorem{lemma}{Lemma}[section]
        \newtheorem{proposition}[lemma]{Proposition}
        \newtheorem{theorem}[lemma]{Theorem}
        \newtheorem{definition}{Definition}[section]
        \newtheorem{remark}[lemma]{Remark}

\title{Stable determination of an inclusion by
boundary measurements\thanks{Work supported by Miur, grant n.2002013279.}}

\author{\begin{tabular}{ccc}
G. Alessandrini\thanks{E-mail: \texttt{alessang@univ.trieste.it}}
& and & M. Di Cristo\thanks{E-mail: \texttt{dicristo@dsm.univ.trieste.it}}
\end{tabular}\\
\normalsize{Dipartimento di Scienze Matematiche}\\
\normalsize{Universit\`a degli Studi di Trieste, Italy}}

\date{}

\begin{document}
\maketitle

\begin{abstract}

We deal with the problem of determining an inclusion within an electrical
conductor from electrical boundary measurements. Under mild a priori assumptions
we establish an optimal stability estimate.
\end{abstract}

\section{Introduction}
In this paper we deal with an inverse boundary value problem which is a special
instance of the well-known Calder\'on's inverse conductivity problem \cite{c}.
Given a bounded domain $\Omega$ in $\mathbb{R}^{n}$, $n\geq2$, with reasonably
smooth boundary, an open set $D$, compactly contained in $\Omega$, and
a constant $k>0$, $k\neq1$, consider, for any $f\in H^{1/2}(\partial\Omega)$,
the weak solution $u\in H^{1}(\Omega)$ to the Dirichlet problem
\begin{eqnarray}
\label{eq}
\textrm{div}((1+(k-1)\chi_{D})\nabla u)&=&0\qquad\textrm{in }\Omega,
\\
\label{dat}
u&=&f\qquad\textrm{on }\partial\Omega,
\end{eqnarray}
where $\chi_{D}$ denotes the characteristic function of the set
$D$. We will denote by $\Lambda_{D}: H^{1/2}(\partial\Omega)\to
H^{-1/2}(\partial\Omega)$ the so called Dirichlet-to-Neumann map,
that is the operator which maps the Dirichlet data onto the
corresponding Neumann data $\frac{\partial
u}{\partial\nu}_{|\partial\Omega}$. The inverse problem that we
examine here is to determine $D$ when $\Lambda_{D}$ is given.

In '88 Isakov \cite{i1} proved the uniqueness, the purpose of the present paper
is to prove a result of stability.
In fact we prove that, under mild a priori assumptions on the regularity
and on the topology of $D$, there is a continuous dependence of $D$
(in the Hausdorff metric) from $\Lambda_{D}$ with a modulus of continuity
of logarithmic type, see Theorem \ref{tp} below.
Let us stress that, indeed, this rate of continuity is the optimal
one, as it was shown by examples in the recent paper \cite{dcr}
by the second author and Luca Rondi.

We wish to mention here
a closely related, but different, problem which attracted a lot
of attention starting from the papers of Friedman \cite{f}
and Friedman and Gustafsson \cite{fg}.
That is the one of determining $D$ when, instead of full knowledge
of the Dirichlet-to-Neumann map, only one, or few, pairs of
Dirichlet and Neumann data are available, see \cite{ai}, \cite{i2}
for extended bibliographical accounts.
Unfortunately, for such a problem, the uniqueness question, not to
mention stability, remains a largely open issue.

Let us illustrate briefly the main steps of our arguments. We must
recall that Isakov's approach to uniqueness is essentially based
on two arguments
\begin{enumerate}
\item[a)]the Runge approximation theorem,
\item[b)]the use of solutions with Green's function type singularities.
\end{enumerate}
Also here we shall use singular solutions, and indeed we shall
need an accurate study of their asymptotic behavior when the
singularity gets close to the set of discontinuity $\partial D$ of
the conductivity coefficient $1+(k-1)\chi_{D}$ in (\ref{eq}), see
Proposition \ref{locbh}. On the other hand, it seems that Runge's
theorem, which is typically based on nonconstructive arguments,
(Lax, \cite{l}, Kohn and Vogelius \cite{kv}) is not suited for
stability estimates and therefore we introduced a different
approach based on quantitative estimates of unique continuation,
see Proposition \ref{p1}.

In Section 2 we formulate our main hypotheses and state the stability
result, Theorem \ref{tp}. In Section 3 we prove Theorem \ref{tp}
on the basis of some auxiliary Propositions, whose proof is deferred to
the following Section 4.

\section{The main result}
\setcounter{equation}{0}
Let us introduce our regularity and topological assumptions on the
conductor $\Omega$ and on the unknown inclusion $D$.
To this purpose we shall need the following definitions.
In places, we shall denote a point $x\in\mathbb{R}^{n}$ by
$x=(x',x_{n})$ where $x'\in\mathbb{R}^{n-1}$, $x_{n}\in\mathbb{R}$.
\begin{definition}
\label{defbound} Let $\Omega$ be a bounded domain in
$\mathbb{R}^{n}$. Given $\alpha$, $0<\alpha\leq1$, we shall say
that a portion $S$ of $\partial\Omega$ is of class $C^{1,\alpha}$
with constants $\overline{r}$, $L>0$ if, for any $P\in S$, there
exists a rigid transformation of coordinates under which we have
$P=0$ and
$$\Omega\cap B_{\overline{r}}(0)=\{x\in B_{\overline{r}}\,:\,
x_{n}>\varphi(x')\},$$
where $\varphi$ is a $C^{1,\alpha}$ function on $B_{\overline{r}}(0)
\subset\mathbb{R}^{n-1}$ satisfying $\varphi(0)=|\nabla\varphi(0)|=0$
and $\Vert \varphi\Vert_{C^{1,\alpha}(B_{\overline{r}}(0))}\leq L\overline{r}$.
\end{definition}
\begin{definition}
We shall say that a portion $S$ of $\partial\Omega$ is of
Lipschitz class with constants $\overline{r}$, $L>0$ if for any
$P\in S$, there exists a rigid transformation of coordinates under
which we have $P=0$ and
$$\Omega\cap B_{\overline{r}}(0)=\{x\in B_{\overline{r}}\,:\,
x_{n}>\varphi(x')\},$$
where $\varphi$ is a Lipschitz continuous function
on $B_{\overline{r}}(0)\subset\mathbb{R}^{n-1}$ satisfying
$\varphi(0)=0$ and  $\Vert \varphi\Vert_{C^{0,1}(B_{\overline{r}}(0))}\leq L\overline{r}$.
\end{definition}
\begin{remark}
We have chosen to scale all norms in a such a way that they are dimensionally equivalent
to their argument. For instance,
for any $\varphi\in C^{1,\alpha}(B_{\overline{r}}(0))$ we set
$$\Vert\varphi \Vert_{C^{1,\alpha}(B_{\overline{r}}(0))}=
\Vert\varphi \Vert_{L^{\infty}(B_{\overline{r}}(0))}
+\overline{r}\Vert\nabla\varphi \Vert_{L^{\infty}(B_{\overline{r}}(0))}
+\overline{r}^{1+\alpha}\vert\nabla\varphi \vert_{\alpha,B_{\overline{r}}(0)}.$$
\end{remark}
\noindent
For given numbers $\overline{r}$, $M$, $\widetilde{\delta}$, $L>0$, $0<\alpha<1$,
we shall assume
\begin{enumerate}
\item[\em{(H1)}] the domain $\Omega$ satisfies the following
conditions
\begin{equation}
\label{ip1}
|\Omega|\leq M\overline{r}^{n},
\end{equation}
where $|\cdot|$ denotes the Lebesgue measure of $\Omega$,
\begin{equation}
\label{ip2}
\partial\Omega\textrm{ is of class }C^{1,\alpha}
\textrm{ with constants }\overline{r},\,L,
\end{equation}
\item[\em{(H2)}] the inclusion $D$ satisfies the following
conditions
\begin{equation}
\label{ip3}
\Omega\smallsetminus\overline{D}\quad\textrm{is connected},
\end{equation}
\begin{equation}
\label{ip4}
\textrm{dist}(D,\partial\Omega)\geq\widetilde{\delta},
\end{equation}
\begin{equation}
\label{ip5}
\partial D\textrm{ is of class }C^{1,\alpha}
\textrm{ with constants }\overline{r},\,L.
\end{equation}
\end{enumerate}
In the sequel we shall refer to numbers $k$, $n$, $\overline{r}$,
$M$, $\widetilde{\delta}$, $L$, $\alpha$
as to the a priori data.
We shall denote by $D_{1}$ and $D_{2}$ two possible inclusions in
$\Omega$, both satisfying the properties mentioned.
We shall denote by $\Lambda_{D_{i}}$, $i=1,2$, the Dirichlet-to-Neumann map
$\Lambda_{D}$ when $D=D_{i}$.
We can now state the main theorem.
\begin{theorem}
\label{tp}
Let $\Omega\subset\mathbb{R}^{n}$, $n\geq2$, satisfy (H1).
Let $k>0$, $k\neq1$ be given.
Let $D_{1}$ and $D_{2}$ be two inclusions in $\Omega$ satisfying (H2).
If, given $\varepsilon>0$, we have
\begin{equation}
\label{b4}
\Vert\Lambda_{D_{1}}-\Lambda_{D_{2}}\Vert_{\mathcal{L}(H^{1/2},H^{-1/2})}
\leq\varepsilon,
\end{equation}
then
$$d_{\mathcal{H}}(\partial{D}_{1},\partial{D}_{2})\leq\omega(\varepsilon),$$
where $\omega$ is an increasing function on $[0,+\infty)$, which satisfies
$$\omega(t)\leq C|\log t|^{-\eta},\qquad\qquad \textrm{for every}\quad 0<t<1$$
and $C$, $\eta$, $C>0$, $0<\eta\leq1$, are constants only depending on the
a priori data.
\end{theorem}
\noindent Here $d_{\mathcal{H}}$ denotes the Hausdorff distance between
bounded closed sets of $\mathbb{R}^{n}$
and $\|\cdot\|_{\mathcal{L}(H^{1/2}H^{-1/2})}$ denotes the operator norm
on the space of bounded linear operators between $H^{1/2}(\partial\Omega)$
and $H^{-1/2}(\partial\Omega)$.
\begin{remark}
It should be emphasized that in this statement the unknown inclusion
may be disconnected.
\end{remark}
\begin{remark}
Several variations of the above results could be devised with
minor adaptations on the arguments. Just to mention one, an
analogous result would be obtained if the Neumann-to-Dirichlet
maps $N_{D_{i}}$ are available instead of the Dirichlet-to-Neumann
maps $\Lambda_{D_{i}}$.
\end{remark}

\section{Proof of Theorem \ref{tp}}
\setcounter{equation}{0}
Before proving Theorem \ref{tp},  we shall state some auxiliary Propositions,
whose proofs are collected in the next Section 4.
Here and in the sequel we shall denote by $\mathcal{G}$ the connected
component of $\Omega\smallsetminus(D_{1}\cup D_{2})$, whose boundary contains
$\partial\Omega$, $\Omega_{D}=\Omega\smallsetminus\overline{\mathcal{G}}$,
$\Omega_{\overline{r}}=\{x\in\mathcal{C}\Omega:\textrm{dist}
(x,\Omega)\leq\overline{r}\}$ and
$\mathcal{S}_{2\overline{r}}=\{x\in\mathbb{R}^{n}:\overline{r}\leq \textrm{dist}(x,\Omega)\leq2
\overline{r}\}$.\\
We introduce a variation of the Hausdorff distance which we call \emph{modified distance}.
\begin{definition}
We shall call modified distance between $D_{1}$ and $D_{2}$ the number
\begin{equation}
\label{dm}
d_{\mu}(D_{1},D_{2})=
\max\bigg\{\sup_{x\in\partial D_{1}\cap\partial\Omega_{D}}\emph{dist}(x,D_{2}),\,
\sup_{x\in\partial D_{2}\cap\partial\Omega_{D}}\emph{dist}(x,D_{1})\bigg\}.
\end{equation}
\end{definition}
\noindent This notion is an adaptation of the one introduced in \cite{abrv},
which was also called \emph{modified distance}.
In order to distinguish such two notions, we call $d_{\mu}$ the present one,
whereas the one in \cite{abrv} was denoted by $d_{m}$.
On the other hand, we need to stress the common peculiarities:
such modified distances do not satisfy the axioms of a metric and in general
do not dominate the Hausdorff distance (see Section 3 in \cite{abrv}
for related arguments).
The following Proposition provides sufficient conditions under which
$d_{\mu}$ dominates $d_{\mathcal{H}}$. See \cite{abrv} Proposition 3.6
for a related statement.
\begin{proposition}
\label{p3}
Let $\Omega$ be an open set in $\mathbb{R}^{n}$ satisfying (H1).
Let $D_{1}$, $D_{2}$ be two bounded open inclusions of $\Omega$ satisfying (H2).
Then
\begin{equation}
\label{cont}
d_{\mathcal{H}}(\partial{D}_{1},\partial{D}_{2})\leq cd_{\mu}(D_{1},D_{2}),
\end{equation}
where $c$ depends only on the a priori assumptions.
\end{proposition}
\noindent With no loss of generality, we can assume that there exists a point $O$
of $\partial D_{1}\cap\partial\Omega_{D}$,
where the maximum in the definition (\ref{dm}) is attained,
that is
\begin{equation}
\label{dPD}
d_{\mu}=d_{\mu}(D_{1},D_{2})=\textrm{dist}(O,D_{2}).
\end{equation}
As is well-known, the Dirichlet-to-Neumann map $\Lambda_{D}$ associated to problem
(\ref{eq}), (\ref{dat}) is defined by:
\begin{equation}
\label{dualp}
<\Lambda_{D}u,v>=\int_{\Omega}(1+(k-1)\chi_{D})\nabla u\cdot\nabla v,
\end{equation}
for every $u\in H^{1}(\Omega)$ solution to (\ref{eq})
and for every $v\in H^{1}(\Omega)$.
Here $<\cdot,\cdot>$ denotes the dual pairing between $H^{-1/2}(\partial\Omega)$
and $H^{1/2}(\partial\Omega)$.
With a slight abuse of notation we shall write
$$<g,f>=\int_{\partial\Omega}gf\,d\sigma,$$
for any $f\in H^{1/2}(\partial \Omega)$
and $g\in H^{-1/2}(\partial \Omega)$.
Let $\Gamma_{D}(x,y)$ be the fundamental solution
for the operator $\textrm{div}((1+(k-1)\chi_{D})\nabla\cdot)$,
thus
\begin{eqnarray}
\label{b1}
&&\textrm{div}((1+(k-1)\chi_{D})\nabla\Gamma_{D}(\cdot,y))=
-\delta(\cdot-y),
\end{eqnarray}
where $y,w\in\mathbb{R}^{n}$, $\delta$ denotes the Dirac distribution .
We shall denote by $\Gamma_{D_{1}}$, $\Gamma_{D_{2}}$
such fundamental solutions when $D=D_{1}$, $D_{2}$ respectively.
Recalling the well-known identity
$$\int_{\Omega}(1+(k-1)\chi_{D_{1}})\nabla u_{1}\cdot\nabla u_{2}-
\int_{\Omega}(1+(k-1)\chi_{D_{2}})\nabla u_{1}\cdot\nabla u_{2}
=\int_{\partial\Omega}u_{1}[\Lambda_{D_{1}}-\Lambda_{D_{2}}]u_{2},$$
which holds
for every $u_{i}\in H^{1}(\Omega)$, $i=1,2$, solutions to (\ref{eq})
when $D=D_{i}$ respectively
(see \cite{i2} formula (5.0.4), Section 5.0), we have
\begin{equation}
\begin{array}{ll}
\label{c1}
&\int_{\Omega}(1+(k-1)\chi_{D_{1}})\nabla\Gamma_{D_{1}}(\cdot,y)\cdot\nabla\Gamma_{D_{2}}
(\cdot,w)
\\[2mm]
&-\int_{\Omega}(1+(k-1)\chi_{D_{2}})\nabla\Gamma_{D_{1}}(\cdot,y)\cdot\nabla
\Gamma_{D_{2}}(\cdot,w)
\\[2mm]
=&\int_{\partial\Omega}\Gamma_{D_{1}}(\cdot,y)
[\Lambda_{D_{1}}-\Lambda_{D_{2}}](\Gamma_{D_{2}}(\cdot,w))d\sigma,
\qquad\forall\,y,w\in\mathcal{C}\overline{\Omega}.
\end{array}
\end{equation}
Let us define, for $y,w\in\mathcal{G}\cup\mathcal{C}\Omega$
\begin{eqnarray}
\label{def1}
S_{D_{1}}(y,w)&=&(k-1)\int_{D_{1}}\nabla\Gamma_{D_{1}}(\cdot,y)
\cdot\nabla\Gamma_{D_{2}}(\cdot,w),
\\[2mm]
\label{def2}
S_{D_{2}}(y,w)&=&(k-1)\int_{D_{2}}\nabla\Gamma_{D_{1}}(\cdot,y)
\cdot\nabla\Gamma_{D_{2}}(\cdot,w),
\\[2mm]
\label{def3}
f(y,w)&=&S_{D_{1}}(y,w)-S_{D_{2}}(y,w).
\end{eqnarray}
Thus (\ref{c1}) can be rewritten as
\begin{equation}
\label{ug}
f(y,w)=\int_{\partial\Omega}\Gamma_{D_{1}}(\cdot,y)
[\Lambda_{D_{1}}-\Lambda_{D_{2}}](\Gamma_{D_{2}}(\cdot,w))d\sigma\qquad\forall\,
y,w\in\mathcal{C}\overline{\Omega}.
\end{equation}
From now on we shall consider the dimension $n\geq3$, since the
case $n=2$ can be treated similarly through minor adaptations
regarding the fundamental solutions. Up to a transformation of
coordinates, we can assume that $O$, defined in (\ref{dPD}), is
the origin of the coordinate system. Let $\nu(O)$ be the outer
unit normal vector to $\partial\Omega_{D}$ in the origin $O$. Such
a normal is indeed well-defined since we are assuming that $O$
realizes the modified distance between $D_{1}$ and $D_{2}$,
therefore, in a small neighborhood of $O$, $\partial\Omega_{D}$ is
made of a part of $\partial D_{1}$, which is known to be
$C^{1,\alpha}$. We will rotate the coordinate system in such a way
that $\nu(O)=(0,\dots,0,-1)$. Taking $y=w=h\nu(O)$, with $h>0$, we
want to evaluate $f(y,y)$ and $S_{D_{1}}(y,y)$ in term of $h$, for
$h$ small. Then, evaluating $S_{D_{2}}$ in term of $d_{\mu}$, we
will get the stability estimate for the modified distance and
thus, using Proposition \ref{p3}, for the Hausdorff distance. An
important ingredient for evaluating $f$ and $S_{D_{1}}$ is the
behavior of the fundamental solution. We state now a proposition
that collects all the results on $\Gamma_{D_{i}}$, $i=1,2$, that
we will need throughout the paper. For $x=(x',x_{n})$, where
$x'\in\mathbb{R}^{n-1}$ and $x_{n}\in\mathbb{R}$, we set
$x^{\star}=(x',-x_{n})$. We shall denote with $\chi^{+}$ the
characteristic function of the half-space $\{x_{n}>0\}$ and with
$\Gamma_{+}$ the fundamental solution of the operator
$\textrm{div}((1+(k-1)\chi^{+})\nabla\cdot)$. If $\Gamma$ is the
standard fundamental solution of the Laplace operator, we have
that (see for instance \cite{aip}, Theorem 4)
\begin{equation}
\label{fondsol}
\Gamma_{+}(x,y)=
\left\{\begin{array}{ll}
\frac{1}{k}\Gamma(x,y)+
\frac{k-1}{k(k+1)}\Gamma(x,y^{\star})\quad&\textrm{for }x_{n}>0,y_{n}>0,\\[2mm]
\frac{2}{k+1}\Gamma(x,y)\quad&\textrm{for }x_{n}y_{n}<0,
,\\[2mm]
\Gamma(x,y)-\frac{k-1}{k+1}\Gamma(x,y^{\star})\quad&\textrm{for }x_{n}<0,y_{n}<0.
\end{array}\right.
\end{equation}
The following Proposition holds.
\begin{proposition}
\label{locbh}
Let $D\subset\mathbb{R}^{n}$ be an open set whose boundary is of class $C^{1,\alpha}$,
with constants $\overline{r}$, $L$.
\begin{enumerate}
\item[(i)]
There exists a constant $c_{1}>0$ depending on $k$, $n$, $\alpha$
and $L$ only, such that
\begin{equation}
\label{lsw}
|\nabla_{x}\Gamma_{D}(x,y)|\leq c_{1}|x-y|^{1-n},
\end{equation}
for every $x,y\in \mathbb{R}^{n}$,
\item[(ii)]
There exist constants $c_{2}$, $c_{3}>0$ depending on $k$, $n$, $\alpha$
and $L$ only, such that
\begin{eqnarray}
\label{locbh0}
&&\big|\Gamma_{D}(x,y)-\Gamma_{+}(x,y)\big|\leq
\frac{c_{2}}{\overline{r}^{\alpha}}|x-y|^{2-n+\alpha},\\[2mm]
\label{locbh1}
&&\big|\nabla_{x}\Gamma_{D}(x,y)-\nabla_{x}\Gamma_{+}(x,y)\big|\leq
\frac{c_{3}}{\overline{r}^{\alpha^{2}}}|x-y|^{1-n+\alpha^{2}},
\end{eqnarray}
for every $x\in D\cap B_{r}(O)$, and for every $y=h\nu(O),$ with
$0<r <\overline{r}_{0}$, $0<h<\overline{r}_{0}$, where
$\overline{r}_{0}=\big(\min\big\{\frac{1}{2}(8L)^{-1/\alpha},\frac{1}{2}\big\}\big)\frac{\overline{r}}{2}
$.
\end{enumerate}
\end{proposition}
\noindent The next two Propositions give us quantitative estimates on $f$ and $S_{D_{1}}$
when we move $y$ towards $O$, along $\nu(O)$.
\begin{proposition}
\label{p1}
Let $\Omega$ be an open set in $\mathbb{R}^{n}$ satisfying (H1).
Let $D_{1},D_{2}$ be two inclusions in $\Omega$ verifying (H2)
and let $y=h\nu(O)$,
with $O$ defined in (\ref{dPD}).
If, given $\varepsilon>0$, we have
$$\Vert\Lambda_{D_{1}}-\Lambda_{D_{2}}\Vert_{\mathcal{L}(H^{1/2},H^{-1/2})}\leq\varepsilon,$$
then for every $h$, $0<h<\overline{c}\,\overline{r}$,
where $0<\overline{c}<1$, depends on $L$,
\begin{equation}
\label{c2}
|f(y,y)|\leq C\frac{\varepsilon^{Bh^{F}}}{h^{A}},
\end{equation}
where $0<A<1$ and $C,B,F>0$ are
constants that depend only on the a priori data.
\end{proposition}
\begin{proposition}
\label{p2}
Let $\Omega$ be an open set in $\mathbb{R}^{n}$ satisfying (H1).
Let $D_{1}$, $D_{2}$ be two
inclusions in $\Omega$ verifying (H2) and $y=h\nu(O)$.
Then for every $h$, $0<h<\overline{r}_{0}/2$,
\begin{equation}
\label{c'3}
|S_{D_{1}}(y,y)|\geq c_{1}h^{2-n}-c_{2}d_{\mu}^{2-2n}+c_{3},
\end{equation}
where $c_{1},c_{2}$ and $c_{3}$ are positive constants only
depending on the a priori data. Here $\overline{r}_{0}$ is the
number introduced in Proposition \ref{locbh}.
\end{proposition}
\noindent Now we have all the tools that we need to
prove Theorem \ref{tp}.
\begin{proof}[Proof of Theorem \ref{tp}.]
Let $O\in \partial D_{1}$ satisfying (\ref{dPD}), that is
$$d_{\mu}(D_{1},D_{2})=\textrm{dist}(O,D_{2})=d_{\mu}.$$
Then, for $y=h\nu(O)$, with $0<h<h_1$,
where $h_1=\min\left\{d_{\mu},\overline{c}\,\overline{r},\overline{r}_{0}/2\right\}$,
using (\ref{lsw}), we have
\begin{equation}
\label{c6}
|S_{D_{2}}(y,y)|\leq c\int_{D_{2}}\frac{1}
{(d_{\mu}-h)^{n-1}}\frac{1}{(d_{\mu}-h)^{n-1}}dx=
c\frac{1}{(d_{\mu}-h)^{2n-2}}|D_{2}|.
\end{equation}
Using Proposition \ref{p1}, we have
\begin{eqnarray}
&&|S_{D_{1}}(y,y)|-|S_{D_{2}}(y,y)|\leq
|S_{D_{1}}(y,y)-S_{D_{2}}(y,y)|
\nonumber\\[2mm]
&=&|f(y,y)|
\leq c\frac{\varepsilon^{Bh^{F}}}{h^{A}}.\nonumber
\end{eqnarray}
On the other hand, by Proposition \ref{p2} and (\ref{c6})
\begin{eqnarray}
|S_{D_{1}}(y,y)|-|S_{D_{2}}(y,y)|
\geq c_{1}h^{2-n}-c_{2}(d_{\mu}-h)^{2-2n}.\nonumber
\end{eqnarray}
Thus we have
$$c_{3}h^{2-n}-c_{4}(d_{\mu}-h)^{2-2n}\leq
\frac{\varepsilon^{Bh^{F}}}{h^{A}}.$$
That is
\begin{eqnarray}
c_{4}(d_{\mu}-h)^{2-2n}&\geq& c_{3}h^{2-n}-\frac{\varepsilon^{Bh^{F}}}{h^{A}}
=h^{2-n}(c_{3}-\varepsilon^{Bh^{F}}h^{\widetilde{A}})
\nonumber\\[2mm]
&\geq&c_{5}h^{2-n}\big(1-\varepsilon^{Bh^{F}}h^{\widetilde{A}}\big),
\end{eqnarray}
where $\widetilde{A}=n-2-A$, $\widetilde{A}>0$. Let
$h=h(\varepsilon)$ where
$h(\varepsilon)=\min\{|\ln\varepsilon|^{-\frac{1}{2F}},
d_{\mu}\}$, for $0<\varepsilon\leq \varepsilon_1$, with
$\varepsilon_{1}$ such that
$\exp(-B|\ln\varepsilon_{1}|^{1/2})=1/2$.
If $d_{\mu}\leq|\ln\varepsilon|^{-\frac{1}{2F}}$
the theorem follows using Proposition \ref{p3}.
In the other case
we have
$$\varepsilon^{Bh(\varepsilon)^{F}}h(\varepsilon)^{\widetilde{A}}
\leq\varepsilon^{B|\ln\varepsilon|^{-1/2}}
\leq\exp\big(-B|\ln\varepsilon|^{1/2}\big).$$
Then, for any $\varepsilon$, $0<\varepsilon<\varepsilon_{1}$,
$$(d_{\mu}-h(\varepsilon))^{2-2n}\geq c_{6}h(\varepsilon)^{2-n},$$
that is
\begin{equation}
\label{c7}
d_{\mu}\leq c_{7}|\ln\varepsilon|^{-\delta\frac{n-2}{2n-2}}
\end{equation}
where $\delta=1/(2F)$. When $\varepsilon\geq\varepsilon_1$, then
$$d_{\mu}\leq\textrm{diam}\,\Omega\leq\textrm{diam}\,\Omega
\frac{|\ln\varepsilon|^{-\frac{1}{2F}}}{|\ln\varepsilon_1|^{-\frac{1}{2F}}}.$$
Finally, using Proposition \ref{p3}, the theorem follows.
\end{proof}

\section{Proofs of the auxiliary Propositions}
\setcounter{equation}{0}
We premise the proof of Proposition \ref{p3} with one lemma.
\begin{lemma}
\label{camm}
Let $\Omega$ be an open set in $\mathbb{R}^{n}$ satisfying (H1).
Let $D$ be a bounded open inclusion of $\Omega$ satisfying (H2).
Then for every $P\in\partial D$, there exists a continuous path $\gamma$
in $\Omega\smallsetminus\overline{D}$ with one end-point in $P$ and
the other on $\partial\Omega$, such that for every $z\in\gamma$
\begin{equation}
\label{dicam}
|z-P|\leq c\,\emph{dist}(z,D),
\end{equation}
where $c$ is a positive constant depending on the a priori data only.
\end{lemma}
\begin{proof}
Using Lemma 5.2 of \cite{abrv},
(which adapted arguments due to Lieberman \cite{li}),
we approximate $\mathrm{dist}(\cdot,\partial D)$
with a regularized distance $\tilde{d}$ such that $\tilde{d}\in C^{2}(\Omega\smallsetminus D)
\cup C^{1,\alpha}(\overline{\Omega\smallsetminus D})$ and the following facts hold
\begin{align*}
&\gamma_{0}\leq\frac{\mathrm{dist}(x,\partial D)}{\tilde{d}(x)}
\leq \gamma_{1},\\
&|\nabla\tilde{d}(y)|\geq c_{1}\quad\textrm{for every }y\in\Omega
\textrm{ s.t. }\mathrm{dist}(y,\partial D)>b\overline{r},\\
&\Vert\tilde{d}\Vert_{1,\alpha}\leq c_{2}\overline{r},
\end{align*}
where $\gamma_{0}$, $\gamma_{1}$, $b$,
$c_{1}$ and $c_{2}$ are positive
constants only depending on $L$ and $\alpha$.
We define for $0<h<a\overline{r}$, with $a$ depending on $L$ and $\alpha$ only,
$$E_{h}=\{x\in\Omega\smallsetminus\overline{D}\,:\,
\tilde{d}(x)>h\}.$$
Arguing as in Lemma 5.3 of \cite{abrv},
$E_{h}$ is connected with boundary of class $C^{1}$ and
\begin{equation}
\label{dset}
\widetilde{c}_{1}h\leq\textrm{dist}(x,\partial{D})
\leq\widetilde{c}_{2}h,\qquad\forall\,x\in\partial E_{h}\cap\Omega,
\end{equation}
where $\widetilde{c}_{1}$, $\widetilde{c}_{2}$ are positive constants
depending on $L$ and $\alpha$ only.
Let us fix $P\in\partial D$. Let $\nu(P)$ be the outer unit normal
to $\partial D$ in $P$ (we recall that $\partial D$ is $C^{1,\alpha}$).
Since (\ref{dset}), there exists a point
$P'\in E_{h}$ such that $P'=\tilde{h}\nu(P)$, where $\tilde{h}$ is a positive constant
$\widetilde{c}_{1}h<\tilde{h}<\widetilde{c}_{2}h$.
We denote by $\overline{PP'}$ the segment whose end-points are $P$ and $P'$.
Since $E_h$ is connected, there exists
a continuous path $\gamma'\subset E_{h}$ with one end-point $P'$ and
the other on $\partial\Omega$.
Since $\gamma'\subset E_{h}$ we have that
for every $x\in\gamma'$, $\textrm{dist}(x,\partial D)\geq ch$, where $c$
is a positive constant.
We then define $\gamma=\gamma'\cup\overline{PP'}$ and the lemma follows.
\end{proof}
\begin{proof}[Proof of Proposition \ref{p3}]
Let us fix $P\in\partial D_{1}$. We distinguish the two following cases.
\begin{enumerate}
\item[i)] $P\in\partial D_{1}\cap\partial\mathcal{G}$,
\item[ii)] $P\in\partial D_{1}\smallsetminus\partial\mathcal{G}$.
\end{enumerate}
If case i) occurs then,
$$\textrm{dist}(P,\partial D_{2})=\textrm{dist}(P,\overline{D}_{2})\leq d_{\mu}.$$
Let us consider case ii). Let $\gamma$ be the continuous path
constructed in Lemma \ref{camm} from $P$ to $\partial\Omega$.
Since $P\notin\partial\mathcal{G}$, there exists
$z\in\gamma\cap\partial D_{2}\cap\partial\Omega_{D}$.
$$\textrm{dist}(z,D_{1})\leq
\sup_{x\in\partial D_{2}\cap\partial\Omega_{D}}\big\{\textrm{dist}(x,D_{1})\big\}
\leq d_{\mu}(D_{1},D_{2}).$$
Thus
$$|z-P|\leq cd_{\mu}(D_{1},D_{2}),$$
where $c>0$ is the constant appearing in (\ref{dicam})
On the other hand
$$\textrm{dist}(P,\partial D_{2})\leq|z-P|.$$
So we obtain that, for every $P\in\partial D_{1}$
$$\mathrm{dist}(P,\partial D_{2})\leq cd_{\mu}(D_{1},D_{2}).$$
Similarly one can show that for every $Q\in\partial D_2$
$$\textrm{dist}(Q,\partial D_{1})\leq cd_{\mu}(D_{1},D_{2}).$$
Then we conclude
$$d_{\mathcal{H}}(\partial D_{1},\partial D_{2})\leq cd_{\mu}(D_{1},D_{2}).$$
\end{proof}

\begin{proof}[Proof of Proposition \ref{locbh}.]
Let us prove (\emph{i}).\\
Let us consider the case $x\in D$ and $y\in\partial D$.
The cases in which $x,y\in D$ or $x,y\in\mathcal{C}D$ are trivial.
Let $h=|x-y|$.
Let $c$ be a positive number less than $\frac{1}{1+2\sqrt{n}}$.
We distinguish the following two cases:
\begin{enumerate}
\item[a)] dist$(x,\partial D)<ch$,
\item[b)] dist$(x,\partial D)\geq ch$.
\end{enumerate}
Let us consider the case a). Let $P\in\partial D$ be such that
$|P-x|=\mathrm{dist}(x,\partial D)$. For every $r>0$, let
$Q_{r}(P)$ be the cube centered at $P$, with sides of length $2r$
and parallel to the coordinates axes. We have that the ball
$B_{r}(P)$ is inscribed into $Q_{r}(P)$. In particular $x\in
Q_{ch}(P)$. On the other hand
$$|P-y|\geq|y-x|-|P-x|\geq h(1-c).$$
Then, due to our choice of $c$, $|P-y|>(2ch)\sqrt{n}$,
that is $y\notin Q_{2ch}(P)$. Thus
$$\textrm{div}_{z}\big((1+(k-1)\chi_{D})\nabla_{z}\Gamma_{D}(z,y)\big)=0
\qquad\textrm{ in }Q_{\frac{3}{2}ch}(P)$$
and for the piecewise $C^{1,\alpha}$ regularity of $\Gamma_{D}$,
proved in \cite{dbef}, see also \cite{lv}, we have
\begin{equation}
\label{1}
\Vert \nabla\Gamma_{D}(\cdot,y)\Vert_{L^{\infty}(Q_{ch}(P))}
\leq\frac{\overline{c}_{1}}{h}\Vert\Gamma_{D}(\cdot,y)\Vert_{L^{\infty}(Q_{\frac{3}{2}ch}(P))},
\end{equation}
where $\overline{c}_{1}$ depends on $L$, $k$, $n$ and $\alpha$
only. Using the pointwise bound of $\Gamma_{D}$ with $\Gamma$ (see
\cite{lsw}), we have
\begin{equation}
\label{lswe}
\Vert\Gamma_{D}(\cdot,y)\Vert_{L^{\infty}(Q_{\frac{3}{2}ch}(P))}\leq\overline{c}_{2}
\left(\frac{ch}{2}\right)^{2-n},
\end{equation}
where $\overline{c}_{2}$ depends on $n$ and $k$ only.
Hence, by (\ref{1}) and (\ref{lswe}), we get
\begin{equation}
\label{2}
|\nabla_{x}\Gamma_{D}(x,y)|
\leq\Vert\nabla\Gamma_{D}(\cdot,y)\Vert_{L^{\infty}(Q_{ch}(P))}
\leq \overline{c}_{3}h^{1-n}
=\overline{c}_{3}|x-y|^{1-n},
\end{equation}
where $\overline{c}_{3}$ depends on $L$, $k$, $n$ and $\alpha$ only.\\
If case b) occurs, then $Q_{\frac{ch}{\sqrt{n}}}(x)\subset D$. Hence
\begin{eqnarray}
&&|\nabla_{x}\Gamma_{D}(x,y)|\leq\Vert\nabla\Gamma_{D}(\cdot,y)
\Vert_{L^{\infty}(Q_{\frac{ch}{2\sqrt{n}}}(x))}
\leq\frac{\overline{c}_{4}}{h}\Vert\Gamma_{D}(\cdot,y)\Vert_{L^{\infty}(Q_{\frac{c}{\sqrt{n}}}(P))}
\nonumber\\
&&\leq\frac{\overline{c}_{4}}{h}(h(1-c))^{2-n}=\overline{c}'_{4}h^{1-n}
=\overline{c}'_{4}|x-y|^{1-n},\nonumber
\end{eqnarray}
where $\overline{c}_{4}$, $\overline{c}'_{4}$
depend on $L$, $k$, $n$ and $\alpha$ only.\\
Let us prove (\emph{ii}).\\
Let us fix $r_{1}=\min\big\{\frac{1}{2}(8L)^{-1/\alpha}\overline{r},\frac{\overline{r}}{2}\big\}$.
Recalling Definition \ref{defbound}, we have that
$$\partial D\cap B_{\overline{r}}(0)=
\{x\in B_{\overline{r}}(0)\,:\,x_{n}=\varphi(x')\},$$
where $\varphi\in C^{1,\alpha}(\mathbb{R}^{n-1})$ satisfying
$\varphi(0)=|\nabla\varphi(0)|=0$.
Let $\theta\in C^{\infty}(\mathbb{R})$ be such that $0\leq\theta\leq1$,
$\theta(t)=1$, for $|t|<1$, $\theta(t)=0$, for $|t|>2$ and
$|\frac{d\theta}{dt}|\leq2$. We consider the following change of
variables $\xi=\Phi(x)$ defined by
$$\left\{\begin{array}{ll}
\xi'=x'\\
\xi_{n}=x_{n}-\varphi(x')\theta\big(\frac{|x'|}{r_{1}}\big)
\theta\big(\frac{x_{n}}{r_{1}}\big).
\end{array}\right.$$
It can be verified that, with the given choice of $r_1$,
the following properties of $\Phi$ hold
\begin{eqnarray}
\label{pr1}
&&\Phi(Q_{2r_{1}}(0))=Q_{2r_{1}}(0),\\[2mm]
\label{pr2}
&&\Phi(Q_{r_{1}}(0)\cap D)=Q_{r_{1}}^{+}(0),\\[2mm]
\label{pr3}
&&c^{-1}|x_{1}-x_{2}|\leq|\Phi(x_{1})-\Phi(x_{2})|\leq c|x_{1}-x_{2}|,
\qquad\forall\,x_{1},x_{2}\in\mathbb{R}^{n},\\[2mm]
\label{pr4}
&&|\Phi(x)-x|\leq\frac{c}{\overline{r}^{\alpha}}|x|^{1+\alpha},
\qquad\forall\,x\in\mathbb{R}^{n},\\[2mm]
\label{pr5}
&&|D\Phi(x)-I|\leq\frac{c}{\overline{r}^{\alpha}}|x|^{\alpha},\qquad\forall\,x\in\mathbb{R}^{n},
\end{eqnarray}
where $Q_{r_{1}}^{+}(0)=\{x\in Q_{r_{1}}(0)\,:\,x_{n}>0\}$
and $c\geq1$ depends on $L$ and $\alpha$ only.
$\Phi$ is a $C^{1,\alpha}$ diffeomorphism from
$\mathbb{R}^{n}$ into itself. Let us define the cylinder $C_{r_{1}}$ as
$$C_{r_{1}}=\{x\in\mathbb{R}^{n}\,:\,|x'|<r_{1},|x_{n}|<r_{1}\}.$$
For $x,y\in C_{r_{1}}$, we have that $\widetilde{\Gamma}_{D}(\xi,\eta)
=\Gamma_{D}(x,y)$, where $\xi=\Phi(x)$, $\eta=\Phi(y)$, is solution of
\begin{equation}
\label{eqnuova}
\textrm{div}_{\xi}((1+(k-1)\chi^{+})B(\xi)\nabla_{\xi}\widetilde{\Gamma}_{D}(\xi,\eta))
=-\delta(\xi-\eta),
\end{equation}
where $B=\frac{JJ^{T}}{\det J}$, with $J=\frac{\partial\xi}{\partial x}(\Phi^{-1}(\xi))$.
We observe that $B$ is of class $C^{\alpha}$ and $B(0)=I$.
Let us consider
$$\widetilde{R}(x,y)=\widetilde{\Gamma}_{D}(x,y)-\Gamma_{+}(x,y),$$
where we keep the notation $x$, $y$ to indicate $\xi$, $\eta$. By
the properties of $\Gamma_+$  and by (\ref{eqnuova}),
$\widetilde{R}$ satisfies
$$\textrm{div}_{x}((1+(k-1)\chi^{+})\nabla_{x}\widetilde{R}(x,y))
=\textrm{div}_{x}((1+(k-1)\chi^{+})(I-B)\nabla_{x}\widetilde{\Gamma}_{D}(x,y)).$$
Let $\widetilde{L}>0$,
depending on the a priori data only,
be such that $\overline{\Omega}\subset B_{\widetilde{L}}(0)$.
Thus using the fundamental solution $\Gamma_{+}$ we obtain
\begin{eqnarray}\nonumber
&&-\widetilde{R}(x,y)=\int_{B_{\widetilde{L}}(0)}(1+(k-1)\chi^{+})(B-I)\nabla_{z}\Gamma_{+}(z,y)
\cdot\nabla_{z}\widetilde{\Gamma}_{D}(z,x)dz\\[2mm]\nonumber
&&+\int_{\partial B_{\widetilde{L}}(0)}(1+(k-1)\chi^{+})\bigg[\widetilde{R}(x,z)
\frac{\partial\Gamma_{+}}{\partial\nu}(z,y)
-\Gamma_{+}(z,y)\frac{\partial \widetilde{R}}{\partial\nu}(x,z)\bigg]d\sigma(z)\\[2mm]\nonumber
&&=\int_{B_{\widetilde{L}}(0)\cap C_{r_{1}}}(1+(k-1)\chi^{+})(B-I)\nabla_{z}\Gamma_{+}(z,y)\cdot
\nabla_{z}\widetilde{\Gamma}_{D}(z,x)dz\\[2mm]\nonumber
&&+\int_{B_{\widetilde{L}}(0)\smallsetminus C_{r_{1}}}(1+(k-1)\chi^{+})
(B-I)\nabla_{z}\Gamma_{+}(z,y)
\cdot\nabla_{z}\widetilde{\Gamma}_{D}(z,x)dz\\[2mm]\nonumber
&&+\int_{\partial B_{\widetilde{L}}(0)}
\bigg[\widetilde{R}(x,z)\frac{\partial\Gamma_{+}}{\partial\nu}(z,y)
-\Gamma_{+}(z,y)\frac{\partial \widetilde{R}}{\partial\nu}(x,z)\bigg]d\sigma(z)\nonumber.
\end{eqnarray}
For $|x|,|y|<r_{1}/2$, the last two integrals are bounded. Using (\ref{lsw})
we obtain
\begin{eqnarray}\nonumber
|\widetilde{R}(x,y)|&\leq& c\bigg(1+\int_{C_{r_{1}}}|z|^{\alpha}
|x-z|^{1-n}|y-z|^{1-n}dz\bigg)\\[2mm]\nonumber
&=&c\bigg(1+I_{1}+I_{2}\bigg),\nonumber
\end{eqnarray}
where $c$ depends on $L$, $\alpha$, $k$ and $n$
and
$$I_{1}=\int_{\{|z|<4h\}\cap C_{r_{1}}}|z|^{\alpha}|x-z|^{1-n}|y-z|^{1-n}dz,$$
$$I_{2}=\int_{\{|z|>4h\}\cap C_{r_{1}}}|z|^{\alpha}|x-z|^{1-n}|y-z|^{1-n}dz.$$
Now
\begin{eqnarray}
I_{1}&\leq&\int_{|w|<4}h^{\alpha}|w|^{\alpha}h^{1-n}\big|\frac{x}{h}-w\big|^{1-n}
h^{1-n}\big|\frac{y}{h}-w\big|^{1-n}h^{n}dw\nonumber\\[2mm]
&=&h^{\alpha+2-n}\int_{|w|<4}|w|^{\alpha}\big|\frac{x}{h}-w\big|^{1-n}
\big|\frac{y}{h}-w\big|^{1-n}dw\nonumber\\[2mm]
&\leq&h^{\alpha+2-n}F(\xi,\eta),\nonumber
\end{eqnarray}
where $h=|x-y|$ and
$$F(\xi,\eta)=4^{\alpha}\int_{|w|<4}|\xi-w|^{1-n}|\eta-w|^{1-n}dw$$
and $\xi=x/h$ and $\eta=y/h$.
From standard bounds (see, for instance, \cite{m}
Chapter 2, Section 11),
it is not difficult to see that
$$F(\xi,\eta)\leq\textrm{const.}<\infty,$$
for all $\xi,\eta\in\mathbb{R}^{n}$, $|\xi-\eta|=1.$
Thus
$$I_{1}\leq c|x-y|^{\alpha+2-n}.$$
Let us consider now $I_{2}$. Since $|y|=-y_{n}\leq|x-y|=h$,
we can deduce $|z|\leq\frac{4}{3}|y-z|$ and $|z|\leq2|x-z|$ and thus
obtain that
$$I_{2}\leq c\int_{|z|>4h}|z|^{\alpha+1-n+1-n}dz\leq ch^{\alpha+2-n}.$$
Then we conclude
\begin{equation}
\label{resto0}
|\widetilde{R}(x,y)|\leq c|x-y|^{\alpha+2-n},
\end{equation}
for every $|x|,|y|<r_{1}/2$, where $c$ depends on $L$, $\alpha$,
$k$ and $n$
only.
Let us go back to the original coordinates system.
We observe that if $x\in\Phi^{-1}(B^{+}_{r_{1}/2}(0))$ and $y=e_{n}y_{n}$,
with $y_{n}\in(-r_{1}/2,0)$ then
$|\Phi(x)-x|$
is bounded by
$c|x-y|^{1+\alpha}$.
Namely, since $\Phi(x)\cdot y\leq 0$ and $\Phi(y)=y$, by (\ref{pr3})
we have
\begin{equation}
\label{coord}
c^{-1}|x|\leq|\Phi(x)|\leq|\Phi(x)-y|\leq c|x-y|.
\end{equation}
On the other hand, by (\ref{pr4}) and (\ref{coord})
\begin{equation}
\label{ve1}
|\Phi(x)-x|\leq\frac{c}{\overline{r}^{\alpha}}|x|^{1+\alpha}
\leq\frac{c'}{\overline{r}^{\alpha}}|x-y|^{1+\alpha}.
\end{equation}
We have
\begin{eqnarray*}
&&R(x,y)=\Gamma_{D}(x,y)-\Gamma_{+}(x,y)\\
&=&\Gamma_{D}(x,y)-\Gamma_{+}(x,y)+\Gamma_{+}(\Phi(x),\Phi(y))-\Gamma_{+}(\Phi(x),\Phi(y))\\
&=&\widetilde{R}(\Phi(x),\Phi(y))+\Gamma_{+}(\Phi(x),y)-\Gamma_{+}(x,y).
\end{eqnarray*}
Using (\ref{pr3}), (\ref{pr4}), (\ref{resto0}) and (\ref{ve1}) we obtain
\begin{eqnarray*}
&&|\Gamma_{D}(x,y)-\Gamma_{+}(x,y)|\\
&\leq&\frac{c}{\overline{r}^{\alpha}}|x-y|^{\alpha+2-n}+
\frac{c}{\overline{r}^{\alpha}}\Vert\nabla\Gamma_{+}(\cdot,y)\Vert_{L^{\infty}(Q_{r_{1}})}
|x-\Phi(x)|\\
&\leq& \frac{c}{\overline{r}^{\alpha}}|x-y|^{\alpha+2-n}+
\frac{c'}{\overline{r}^{\alpha}}|x-y|^{1+\alpha}h^{1-n}\\
&\leq& \frac{c''}{\overline{r}^{\alpha}}|x-y|^{\alpha+2-n},
\end{eqnarray*}
where $c''$ depends on $k$, $n$, $\alpha$ and $L$
only.
We estimate now the first derivative of $R$. To estimate the first derivative
of $\widetilde{R}$ let us consider
a cube $Q\subset B_{r_{1}/4}^{+}(x)$
of side $cr_{1}/4$, with $0<c<1$,
such that $x\in\partial Q$.
The following interpolation inequality holds:
$$\Vert\nabla \widetilde{R}(\cdot,y)\Vert_{L^{\infty}(Q)}\leq
c\Vert \widetilde{R}(\cdot,y)\Vert^{1-\delta}_{L^{\infty}(Q)}
|\nabla \widetilde{R}(\cdot,y)|_{\alpha,Q}^{\delta},$$
where $\delta=\frac{1}{1+\alpha}$,
$c$ depends on $L$ only and
$$|\nabla \widetilde{R}|_{\alpha,Q}=\sup_{x,x'\in Q,x\neq x'}
\frac{|\nabla \widetilde{R}(x,y)-\nabla
\widetilde{R}(x',y)|}{|x-x'|^{\alpha}}.$$ Since, from the
piecewise H\"older continuity of $\nabla\Gamma_{D}$ see (\ref{1}),
and also of $\nabla\Gamma_{+}$, see (\ref{fondsol}), we have that
$$|\nabla \widetilde{R}(\cdot,y)|_{\alpha,Q}\leq
|\nabla\widetilde{\Gamma}_{D}(\cdot,y)|_{\alpha,Q}+|\nabla \Gamma_{+}(\cdot,y)|_{\alpha,Q}
\leq ch^{-\alpha+1-n},$$
where $c$ depends on $L$ only, thus
we conclude
$$|\nabla_{x}\widetilde{R}(x,y)|\leq \frac{c}{\overline{r}^{\eta}}
h^{(\alpha+2-n)(1-\delta)}h^{(-\alpha+1-n)\delta}
=\frac{c}{\overline{r}^{\eta}} h^{1-n+\eta},
$$
where $\eta=\frac{\alpha^2}{1+\alpha}$. Thus
\begin{equation}
\label{resto1}
|\nabla_{x}\widetilde{R}(x,y)|\leq \frac{c}{\overline{r}^{\eta}}|x-y|^{\eta+1-n},
\end{equation}
where $\eta=\frac{\alpha^{2}}{1+\alpha}$ and
$c$ depends on $L$ only.
Concerning $\Gamma_{+}$ we have
\begin{eqnarray*}
&&|\nabla_{x}\Gamma_{+}(\Phi(x),y)-\nabla_{x}\Gamma_{+}(x,y)|\\[2mm]
&=&|D\Phi(x)^{T}\nabla\Gamma_{+}(\cdot,y)_{|\Phi(x)}-\nabla_{x}\Gamma_{+}(x,y)|\\[2mm]
&\leq&|(D\Phi(x)^{T}-I)\nabla\Gamma_{+}(\cdot,y)_{|\Phi(x)}|\\
&&+|\nabla\Gamma_{+}(\cdot,y)_{|\Phi(x)}-\nabla_{x}\Gamma_{+}(x,y)|\\
&\leq& \frac{c}{\overline{r}^{\alpha}}\Vert\nabla\Gamma_{+}(\cdot,y)\Vert_{L^{\infty}(Q_{r_{1}})}
|x-\Phi(x)|+
|\nabla\Gamma_{+}(\cdot,y)|_{\alpha,Q}|\Phi(x)-x|^{\alpha}\\
&\leq&\frac{c'}{\overline{r}^{\alpha}}h^{1+\alpha}h^{1-n}+
 \frac{c}{\overline{r}^{\alpha^{2}}}h^{-\alpha+1-n}h^{(1+\alpha)\alpha}
\\[2mm]
&\leq&\frac{c}{\overline{r}^{\alpha^{2}}}h^{1-n+\alpha^{2}},
\end{eqnarray*}
where $c$ depends on $k$, $n$, $\alpha$ and $L$
only.
\end{proof}

\begin{proof}[Proof of Proposition \ref{p1}.]
Let us fix $\overline{y}\in\mathcal{S}_{2\overline{r}}$ and let us consider
$f(\overline{y},\cdot)$. We have that
\begin{equation}
\label{d1}
\Delta_{w}f(\overline{y},w)=0\qquad\textrm{ in }\mathcal{C}\overline{\Omega}_{D}.
\end{equation}
For $w\in\mathcal{S}_{2\overline{r}}$, by (\ref{b4}), (\ref{ug}) and (\ref{lsw}) we have
\begin{equation}
\label{d2}
|f(\overline{y},w)|\leq C(\overline{r},L,M)\Vert\Lambda_{D_{1}}-\Lambda_{D_{2}}\Vert=
\widetilde{\varepsilon}.
\end{equation}
Let us now estimate $f(\overline{y},w)$ when $w\in\mathcal{G}$. We define
$\mathcal{G}^{h}=\{x\in\mathcal{G}:\textrm{dist}(x,\Omega_{D})\geq h\}$. For every
$w\in\mathcal{G}^{h}$, we have that
\begin{eqnarray}
\label{d3}
|S_{D_{1}}(\overline{y},w)|&\leq&|k-1|\int_{D_{1}}|\nabla_{x}\Gamma_{D_{1}}(x,\overline{y})|\,
|\nabla_{x}\Gamma_{D_{2}}(x,w)|dx
\nonumber\\[2mm]
&\leq& c\int_{D_{1}}|x-w|^{1-n}dx\leq ch^{1-n}.
\end{eqnarray}
Similarly $|S_{D_{2}}(\overline{y},w)|\leq ch^{1-n}$. Then we conclude that
\begin{equation}
\label{d4}
|f(\overline{y},w)|\leq ch^{1-n}
\qquad\textrm{ in }\mathcal{G}^{h}.\nonumber
\end{equation}
At this stage we shall make use of the three spheres inequality
for supremum norms of harmonic functions $v$, see for instance
\cite{km}, \cite{k}. For every $l_{1}$, $l_{2}$, $1<l_{1}<l_{2}$
and for every
$x\in\mathcal{G}\cup\mathcal{S}_{2\overline{r}}\cup\Omega_{\overline{r}}$
there exists $\tau\in(0,1]$, depending only on $l_{1}$, $l_{2}$
and $n$ such that
$$\Vert v\Vert_{L^{\infty}(B_{l_{1}r}(x))}
\leq\Vert v\Vert^{\tau}_{L^{\infty}(B_{r}(x))}
\Vert v\Vert^{1-\tau}_{L^{\infty}(B_{l_{2}r}(x))}.$$
We apply it for $v(\cdot)=f(\overline{y},\cdot)$ in the ball $B_{\overline{r}}(\overline{x})$,
where $\overline{x}\in\mathcal{S}_{2\overline{r}}$
be such that $\textrm{dist}(\overline{x},\Gamma)=\overline{r}/2$,
where $\Gamma=\{x\in\mathbb{R}^{n}:\textrm{dist}(x,\Omega)=
\overline{r}\}\subset\partial\mathcal{S}_{2\overline{r}}$, $l_{1}=3r=3\overline{r}/2$ and
$l_{2}=4r=2\overline{r}$, then we obtain
\begin{equation}
\label{d5}
\Vert f(\overline{y},\cdot)\Vert_{L^{\infty}(B_{3\overline{r}/2}(\overline{x}))}
\leq\Vert f(\overline{y},\cdot)\Vert^{\tau}_{L^{\infty}(B_{\overline{r}/2}(\overline{x}))}
\Vert f(\overline{y},\cdot)\Vert^{1-\tau}_{L^{\infty}(B_{2\overline{r}}(\overline{x}))}.
\end{equation}
For every $\overline{w}\in\mathcal{G}^{h}$, we denote with $\gamma$ a
simple arc in $\overline{\mathcal{G}}\cup
\overline{\Omega}_{\overline{r}}\cup\overline{\mathcal{S}}_{2\overline{r}}$
joining $\overline{x}$ to $\overline{w}$.
Let us define $\{x_{i}\}$, $i=1,\dots,s$ as follows
$x_{1}=\overline{x}$, $x_{i+1}=\gamma(t_{i})$, where
$t_{i}=\max\{t:|\gamma(t)-x_{i}|=\overline{r}\}$ if
$|x_{i}-\overline{w}|>\overline{r}$, otherwise let $i=s$ and stop the
process. By construction, the balls $B_{\overline{r}/2}(x_{i})$ are
pairwise disjoint, $|x_{i+1}-x_{i}|=\overline{r}$ for $i=1,\dots,s-1$,
$|x_{s}-\overline{w}|\leq\overline{r}$. For (\ref{ip1}), there exists $\beta$
such that $s\leq\beta$. An iterated
application of the three spheres inequality (\ref{d5}) for
$f(\overline{y},\cdot)$ (see for instance \cite{abrv} pg.780,
\cite{adb} Appendix E) gives that
for any $r$, $0<r<\overline{r}$
\begin{equation}
\label{d6}
\Vert f(\overline{y},\cdot)\Vert_{L^{\infty}(B_{r/2}(\overline{w}))}
\leq\Vert f(\overline{y},\cdot)\Vert^{\tau}_{L^{\infty}(B_{r/2}(\overline{x}))}
\Vert f(\overline{y},\cdot)\Vert^{1-\tau}_{L^{\infty}(\mathcal{G})}.
\end{equation}
We can now estimate the right hand side of (\ref{d6}) by (\ref{d2}) and
(\ref{d3}) and obtain, for any $r$, $0<r<\overline{r}$
\begin{equation}
\label{d7}
\Vert f(\overline{y},\cdot)\Vert_{L^{\infty}(B_{r/2}(\overline{w}))}
\leq
c(h^{1-n})^{1-\tau^{s}}\varepsilon^{\tau^{s}}
\leq c(h^{1-n})^{A}\varepsilon^{\widetilde{\beta}},
\end{equation}
where $\widetilde{\beta}=\tau^{\beta}$ and $A=1-\widetilde{\beta}$.
Let $O\in\partial D_{1}$, as defined in (\ref{dPD}), that is
$$d(O,D_{2})=d_{\mu}(D_{1},D_{2}).$$
There exists a $C^{1,\alpha}$ neighborhood $U$ of $O$ in
$\partial\Omega_{D}$ with constants $\overline{r}$ and $L$.
Thus there exists a non-tangential vector field
$\widetilde{\nu}$, defined on $U$
such that the truncated cone
\begin{equation}
\label{cone}
C(O,\widetilde{\nu}(O),\theta,\overline{r})=\bigg\{x\in\mathbb{R}^{n}:
\frac{(x-O)\cdot\widetilde{\nu}(O)}{|x-O|}>\cos\theta,\,|x-O|<\overline{r}\bigg\}
\end{equation}
satisfies
$$C(O,\widetilde{\nu}(O),\theta,\overline{r})\subset\mathcal{G},$$
where $\theta=\arctan(1/\overline{L})$. Let us define
\begin{eqnarray}
&&\lambda_{1}=\min\bigg\{\frac{\overline{r}}{1+\sin\theta},
\frac{\overline{r}}{3\sin\theta}\bigg\},
\nonumber\\[2mm]
&&\theta_{1}=\arcsin\bigg(\frac{\sin\theta}{4}\bigg),
\nonumber\\[2mm]
&&w_{1}=O+\lambda_{1}\nu,
\nonumber\\[2mm]
&&\rho_{1}=\lambda_{1}\sin\theta_{1}.\nonumber
\end{eqnarray}
We have that $B_{\rho_{1}}(w_{1})\subset C(O,\widetilde{\nu}(O),\theta_{1},\overline{r})$,
$B_{4\rho_{1}}(w_{1})\subset C(O,\widetilde{\nu}(O),\theta,\overline{r})$. Let
$\overline{w}=w_{1}$, since $\rho_{1}\leq\overline{r}/2$, we can use
(\ref{d7}) in the ball $B_{\rho_{1}}(\overline{w})$ and we can approach
$O\in\partial D_{1}$ by constructing a sequence of balls contained in the cone
$C(O,\widetilde{\nu}(O),\theta_{1},\overline{r})$. We define, for $k\geq2$
$$w_{k}=O+\lambda_{k}\nu,\qquad
\lambda_{k}=\chi\lambda_{k-1},\qquad
\rho_{k}=\chi\rho_{k-1},\quad\textrm{ with }
\chi=\frac{1-\sin\theta_{1}}{1+\sin\theta_{1}}.$$
Hence $\rho_{k}=\chi^{k-1}\rho_{1}$, $\lambda_{k}=\chi^{k-1}\lambda_{1}$ and
$$B_{\rho_{k+1}}(w_{k+1})\subset B_{\rho_{3k}}(w_{k})\subset
B_{\rho_{4k}}(w_{k})\subset C(O,\nu,\theta,\overline{r}).$$
Denoting $d(k)=|w_{k}-O|-\rho_{k}=\lambda_{k}-\rho_{k}$, we have
$d(k)=\chi^{k-1}d(1)$, with $d(1)=\lambda_{1}(1-\sin\theta)$.
For any $r$, $0<r\leq d(1)$, let $k(r)$ be the smallest integer such that
$d(k)\leq r$, that is
$$\frac{\big|\log\frac{r}{d(1)}\big|}{\big|\log\chi\big|}
\leq k(r)-1\leq\frac{\big|\log\frac{r}{d(1)}\big|}{\big|\log\chi\big|}+1.$$
By an iterated application of the three spheres inequality over the
chain of balls $B_{\rho_{1}}(w_{1}),\dots,B_{\rho_{k(r)}}(w_{k(r)})$,
we have
\begin{eqnarray}
\label{d9}
\Vert f(\overline{y},\cdot)\Vert_{L^{\infty}(B_{\rho_{k(r)}}(w_{k(r)}))}
&\leq& c(h^{1-n})^{A(1-\tau^{k(r)-1})}
\varepsilon^{\beta\tau^{k(r)-1}}
\nonumber\\
&\leq&c(h^{1-n})^{A}\varepsilon^{\beta\tau^{k(r)-1}},
\end{eqnarray}
for $0<r<c\overline{r}$, where $0<c<1$ depends on $L$ only.
$$ $$
Let us consider now $f(y,w)$ as a function of $y$.
First observe that
\begin{equation}
\label{d11}
\Delta_{y}f(y,w)=0\qquad\textrm{in }\mathcal{C}\Omega_{D},
\qquad\textrm{for all }w\in\mathcal{C}\Omega_{D}.\nonumber
\end{equation}
For $y,w\in\mathcal{G}^{h}$, $y\neq w$,
using (\ref{lsw}),
we have
$$|S_{D_{1}}(y,w)|\leq c\int_{D_{1}}|x-y|^{1-n}|x-w|^{1-n}dx
\leq ch^{2-n}.$$
Similarly for $S_{D_{2}}$. Therefore
\begin{equation}
\label{d12}
|f(y,w)|\leq ch^{2-2n}\qquad\textrm{ with }y,w\in\mathcal{G}^{h}.\nonumber
\end{equation}
Finally, for $y\in\mathcal{S}_{2\overline{r}}$ and
$w\in\mathcal{G}^{h}$,
using (\ref{d9}), we have
\begin{equation}
\label{d13}
|f(y,w)|\leq c(h^{1-n})^{A}\varepsilon^{\beta\tau^{k(h)-1}}.\nonumber
\end{equation}
Proceeding as before, let us fix $w\in\mathcal{G}$ such that
$\textrm{dist}(w,\partial\Omega_{D})=h$ and
$\widetilde{y}\in\mathcal{S}_{2\overline{r}}$ such that
$\textrm{dist}(\widetilde{y},\Gamma)=\overline{r}/2$. Taking
$r=\overline{r}/2$, $l_{1}=3r$, $l_{2}=4r$,
$y_{1}=O+\lambda_{1}\nu$ and using iteratively the three spheres
inequality, we have
$$\Vert f(y,w)\Vert_{L^{\infty}(B_{\overline{r}/2}(y_{1}))}
\leq \Vert f(y,w)\Vert^{\tau^{s}}_{L^{\infty}(B_{\overline{r}/2}(\widetilde{y}))}
\Vert f(y,w)\Vert^{1-\tau^{s}}_{L^{\infty}(\mathcal{G})},$$
where $\tau$ and $s$ are the same number established previously.
Therefore
\begin{eqnarray}
\label{d14}
\Vert f(y,w)\Vert^{\tau^{s}}_{L^{\infty}(B_{\overline{r}/2}(y_{1})))}
&\leq&
c(h^{2-2n})^{1-\tau^{s}}(h^{1-n})^{A\tau^{s}}
(\varepsilon^{\beta\tau^{k(h)-1}})^{\tau^{s}}
\nonumber\\[2mm]
&\leq&c(h^{2-2n})^{1-\gamma}(h^{1-n})^{A\tau^{s}}
(\varepsilon^{\beta\tau^{k(h)-1}}
)^{\gamma}
\nonumber\\[2mm]
&\leq& c(h^{2-2n})^{A'}(\varepsilon^{\beta\tau^{k(h)-1}}
)^{\gamma},\nonumber
\end{eqnarray}
where $\gamma=\tau^{\beta}$, with $\beta$ as before, so
$0<\gamma<1$, and $A'=A\tau^{s}+1-\gamma$. Once more, let us apply
iteratively the three spheres inequality over a chain of balls
contained in a cone with vertex in $O$ and we obtain
\begin{equation}
\label{d15}
\Vert f(y,w)\Vert_{L^{\infty}(B_{\rho_{k}}(y_{k(h)})}
\leq
c(h^{2-2n})^{A'(1-\tau^{k(h)-1})}(\varepsilon^{\beta\tau^{k(h)-1}}
)^{\gamma\tau^{k(h)-1}}.
\end{equation}
Now, from (\ref{d15}), choosing $y=w=h\nu(O)$,
where $\nu(O)$ is the exterior unit normal to $\partial\Omega_{D}$
in $O$, we obtain
\begin{eqnarray}
\label{holderf}
|f(y,y)|
\leq ch^{A''}(\varepsilon^{\beta\tau^{k(h)-1}}
)^{\gamma\tau^{k(h)-1}},
\end{eqnarray}
where $A''=-(2-2n)\beta A'>0$. We observe that,
for $0<h<c\overline{r}$,
where $0<c<1$ depends on $L$,
$k(h)\leq c|\log h|=-c\log h$, so we can write
$$\tau^{k(h)}=\textrm{e}^{-c\log h\log(\tau)}
=h^{-c\log\tau}=h^{c|\log\tau|}=h^{F},$$
with $F=c|\log\tau|$.
Therefore
\begin{eqnarray}
|f(y,y)|&\leq& h^{-A''}\varepsilon^{B\tau^{k(h)}}
\nonumber\\[2mm]
&=&\textrm{e}^{-A''\log h}\textrm{e}^{B\tau^{k(h)}\log\varepsilon}
\nonumber\\[2mm]
&=&\textrm{e}^{-A''\log h+B'h^{F}\log\varepsilon}.\nonumber
\end{eqnarray}
Then in (\ref{holderf}) we obtain that
\begin{eqnarray}
|f(y,y)|&\leq&\textrm{e}^{-A'\log h+B' h^{F}\log\varepsilon}
=\frac{\varepsilon^{B'h^{F}}}{h^{A'}}.\nonumber
\end{eqnarray}
\end{proof}
\begin{proof}[Proof of Proposition \ref{p2}.]
Let us consider $y=h\nu(O)$, where $\nu(O)$ is the
exterior outer normal to $\partial\Omega_{D}$ in $O$
with $O$ defined as in (\ref{dPD}),
$0<h<
\overline{r}_{0}$, where $\overline{r}_0$ is the number introduced
in Proposition \ref{locbh} and $x\in D_{1}$ such that $|x-y|<r$,
with $0<r<\overline{r}_{0}$. Let us first observe that since
$O\in\partial D_{1}$ and $x\in D_{1}$, for $\Gamma_{D_{1}}$ we
have the asymptotic formula (\ref{locbh1}), which says that
$$\bigg|\nabla_{x}\Gamma_{D_{1}}(x,y)
-\nabla_{x}\Gamma_{+}(x,y)\bigg|\leq c_{1}|x-y|^{1-n+\delta}.$$
Furthermore, since we are in the situation in which $x\in D_{1}$
and $y\notin D_{1}$, for (\ref{fondsol}),
$\Gamma_{+}(x,y)=2/(k+1)\Gamma(x,y)$, where $\Gamma(x,y)$ denotes
the standard fundamental solution of the Laplace operator. Let us
consider now $\Gamma_{D_{2}}(x,y)$. With our choice of $O$, $x$
and $y$, we know that $y\notin D_{2}$ but we do not have any
information on $x$, that is we do not know in which side of the
interface $\partial D_{2}$ it is.
Thus we have to distinguish different situations.\\
If $x\in B_{r}(O)\cap D_{1}\cap D_{2}$, then we have the asymptotic
formula (\ref{fondsol}) for $\Gamma_{D_{2}}$ and
from Lemma 3.1 of \cite{a2} the following formula holds
\begin{equation}
\label{sottostima}
\nabla_{x}\Gamma_{D_{1}}(x,y)\cdot\nabla_{x}\Gamma_{D_{2}}(x,y)
\geq c|x-y|^{2-2n}.
\end{equation}
Consider now the case $x\in(D_{1}\smallsetminus D_{2})\cap B_{r}(O)$.
In this region let us consider a smaller ball $B_{\rho}(O)$
centered in $O$ with radius $\rho$ where $0<\rho< d_{\mu}$.
Since the definition of $d_{\mu}$ we have $B_{\rho}\cap D_{2}=\emptyset$.
If $x$ and $y$ are in $B_{\rho}(O)$, we have
\begin{equation}
\label{max}
\left\{\begin{array}{l}
\Delta_{x}\big(\Gamma_{D_{2}}(x,y)-\Gamma(x,y)\big)=0\qquad\textrm{in }B_{\rho}(O),\\[2mm]
\big[\Gamma_{D_{2}}(x,y)-\Gamma(x,y)\big]_{|\partial B_{\rho}(O)}\leq c\rho^{2-n}.
\end{array}\right.
\end{equation}
Thus by the maximum principle
\begin{equation}
\label{max0}
\big|\Gamma_{D_{2}}(x,y)-\Gamma(x,y)\big|\leq c_{1}\rho^{2-n}\qquad\forall\,x,y\in B_{\rho}(O),
\end{equation}
and by interior gradient bound
\begin{equation}
\label{max1}
\big|\nabla_{x}\Gamma_{D_{2}}(x,y)-\nabla_{x}\Gamma(x,y)\big|
\leq c_{2}\rho^{1-n}\qquad\forall\,
x\in B_{\rho/2}(O),\forall\,y\in B_{\rho}(O).
\end{equation}
Thus, using Lemma 3.1 of \cite{a2}, in $B_{\rho/2}(O)$ we obtain the formula
\begin{equation}
\label{sottstima2}
\nabla_{x}\Gamma_{D_{1}}(x,y)\cdot\nabla_{x}\Gamma_{D_{2}}(x,y)\geq
c_{3}|x-y|^{2-2n}-c_{4}\rho^{2-2n}.
\end{equation}
Let us consider
$h\leq \overline{r}_{0}/2$ and $B_{r}(O)
=\{x\in\mathbb{R}^{n}:|x-O|<r\}$, with $0<r<\overline{r}_{0}$. Then we have
\begin{eqnarray}
&&|S_{D_{1}}(y,y)|\nonumber\\[2mm]
&=&|k-1|\bigg|\int\limits_{D_{1}\cap B_{r}(O)}
\nabla\Gamma_{D_{1}}\cdot\nabla\Gamma_{D_{2}}dx+\int\limits_{D_{1}\smallsetminus
B_{r}(O)}\nabla\Gamma_{D_{1}}\cdot\nabla\Gamma_{D_{2}}dx\bigg|
\nonumber\\[2mm]
&\geq&|k-1|\bigg|\int\limits_{D_{1}\cap B_{r}(O)}
\nabla\Gamma_{D_{1}}\cdot\nabla\Gamma_{D_{2}}dx\bigg|-|k-1|\bigg|\int\limits_{D_{1}\smallsetminus
B_{r}(O)}\nabla\Gamma_{D_{1}}\cdot\nabla\Gamma_{D_{2}}dx\bigg|
\nonumber
\end{eqnarray}
The first term can be estimated as follows
\begin{eqnarray}
&&\bigg|\int_{D_{1}\cap B_{r}(O)}
\nabla\Gamma_{D_{1}}\cdot\nabla\Gamma_{D_{2}}dx\bigg|\nonumber\\[2mm]
&=&\bigg|\int_{(D_{1}\cap D_{2})\cap B_{r}(O)}
\nabla\Gamma_{D_{1}}\cdot\nabla\Gamma_{D_{2}}dx+
\int_{(D_{1}\smallsetminus D_{2})\cap B_{r}(O)}
\nabla\Gamma_{D_{1}}\cdot\nabla\Gamma_{D_{2}}dx\bigg|\nonumber\\[2mm]
&\geq&\bigg|\int_{(D_{1}\cap D_{2})\cap B_{r}(O)}
\nabla\Gamma_{D_{1}}\cdot\nabla\Gamma_{D_{2}}dx+
\int_{(D_{1}\smallsetminus D_{2})\cap B_{\rho}(O)}
\nabla\Gamma_{D_{1}}\cdot\nabla\Gamma_{D_{2}}dx\bigg|\nonumber\\[2mm]
&&-\bigg|\int_{[(D_{1}\smallsetminus D_{2})\cap B_{r}(O)]\smallsetminus B_{\rho}(O)}
\nabla\Gamma_{D_{1}}\cdot\nabla\Gamma_{D_{2}}dx
\bigg|\nonumber
\end{eqnarray}
In conclusion, choosing $\rho=d_{\mu}/2$ and
using (\ref{sottostima}) and (\ref{lsw}) we obtain
\begin{eqnarray}
|S_{D}(y,y)|
&\geq&c_{1}\int\limits_{[(D_{1}\cap D_{2})\cap B_{r}(O)]\cup[(D_{1}\smallsetminus
D_{2})\cap B_{d_{\mu}/2}(O)]}|x-y|^{2-2n}dx
\nonumber\\[2mm]
&&
-c_{2}\int\limits_{[(D_{1}\smallsetminus D_{2})\cap B_{r}(O)]
\smallsetminus B_{d_{\mu}/2}(O)}
|x-y|^{1-n}|x-y|^{1-n}dx
\nonumber\\[2mm]
&&-c_{3}\int\limits_{D_{1}\smallsetminus B_{r}(O)}
|x-y|^{1-n}|x-y|^{1-n}dx
\nonumber\\[2mm]
&\geq& c_{4}h^{2-n}-c_{5}d_{\mu}^{2-2n}-c_{7}.\nonumber
\end{eqnarray}
\end{proof}

{\bf Acknowledgements.} The authors wish to express their
gratitude to Professors Edi Rosset and Sergio Vessella for
fruitful discussions on the topics of this paper.


\begin{thebibliography}{} \label{bbiibb}
\bibitem[A]{a2}        G. Alessandrini:
                        {\em Singular solutions of elliptic equations and
                             the determination of conductivity by boundary
                             measurements,}
                        J. Differential Equations, 84, 1990, pp.252-272.
\bibitem[A-B-R-V]{abrv} G. Alessandrini, E. Beretta, E. Rosset, S.Vessella:
                        {\em Optimal stability for inverse elliptic boundary
                             value problems with unknown boundaries.}
                        Ann. Scuola Norm. Sup. Pisa Cl. Sci., 4, XXIX, 2000,
                        pp.755-806.
\bibitem[A-DB]{adb}  G. Alessandrini, E. Di Benedetto:
                     {\em Determining 2-dimensional cracks in 3-dimensional
                          bodies: uniqueness and stability.}
                     Indiana Univ. Math. J., 46, 1997, pp.1-82.
 \bibitem[A-I]{ai}  G. Alessandrini, V. Isakov:
                     {\em Analicity and uniqueness for the inverse
                         conductivity problem.}
                    Rend. Istit. Mat. Univ. Trieste, 28, no.1-2,1997,
pp.351-369.
\bibitem[A-I-P]{aip}  G. Alessandrini, V. Isakov, J. Powell:
                      {\em Local uniqueness in the inverse conductivity
                           problem with one measurement.}
                      Trans. Amer. Math. Soc., 347, 1995, pp.3031-3041.
\bibitem[C]{c}  A.P. Calder\'on:
                      {\em On an inverse boundary value problem.}
                      Seminar on numerical analysis and its applications to continuum
                       physics, Societade Brasileira de Matem\'atica,
                       Rio de Janeiro, 1980, pp.65-73.
\bibitem[DB-E-F]{dbef} E. Di Benedetto, C. Elliott, A. Friedman:
                      {\em The free boundary of a flow in a porous body
                        heated from its boundary.}
                      Nonlinear Anal., 10, no.9, 1986, pp.879-900.
\bibitem[DC-R]{dcr} M. Di Cristo, L. Rondi :
                      {\em Examples of exponential instability for inverse inclusion
                           and scattering problem.}
                      Inverse Problems, 19, no.3, 2003, pp.685-701.
\bibitem[F]{f}  A. Friedman:
                      {\em Detection of mines by electric.}
                      SIAM J. Appl. Math, 47, no.1, 1987, pp.201-212.
\bibitem[F-G]{fg}  A. Friedman, B.Gustafsson:
                      {\em Identification of the conductivity coefficient in an
                        elliptic equation.}
                      SIAM J. Math Anal., 18, no.18, 1987, pp.777-787.
\bibitem[I1]{i1} V. Isakov:
                      {\em On uniqueness of recovery of a discontinuous
                       conductivity coefficient.}
                      Comm. Pure Appl. Math., 41, no.7, 1988, pp.865-877.
\bibitem[I2]{i2}   V. Isakov:
                      {\em Inverse problems for partial differential equations.}
                     Springer Verlag, 1998.
\bibitem[K-V]{kv}  R. Kohn. M. Vogelius
                      {\em Determining conductivity by boundary measurements, II,
                           interior results.}
                      Comm. Pure Appl. Math., 38, 1985, pp.643-667.
\bibitem[K-M]{km}      J. Korevaar, J. Meyers:
                      {\em Logarithmic convexity for supremum norms of harmonic functions.}
                      Bull. London Math. Soc., 26, 1994, pp.353-362.
\bibitem[K]{k}        I. Kukavica:
                      {\em Quantitative uniqueness for second-order elliptic
                           operators.}
                      Duke Math. J., 91, 1998, pp.225-240.
\bibitem[L]{l}        P. Lax
                      {\em A Stability theorem for solutions of abstract differential equations and its
                      applications to the study of the local behaviour of solutions of elliptic equations.}
                      Comm. Pure Appl. Math., 9, 1956, pp.747-766.
\bibitem[L-V]{lv}     Y.Y.Li, M.Vogelius:
                      {\em Gradient estimates for solutions to divergence form elliptic equations
                        with discontinuous coefficients.}
                      Arch. Rational Mech. Anal., 153, 2000, pp.91-151.
\bibitem[Li]{li}       G.M. Lieberman:
                      {\em Regularized distance and its applications.}
            Pacific J. Math., 117, 1985, pp.329-353.
\bibitem[L-S-W]{lsw}  W. Littman, G. Stampacchia, H. Weinberger:
                      {\em Regular points for elliptic equations with discontinuous coefficients.}
                      Ann. Scuola Norm. Sup. Pisa, 4, XXVII, 1963, pp.43-77.
\bibitem[M]{m}        C. Miranda:
                      {\em Partial differential equations of elliptic type.}
                      Springer-Verlag, New York-Berlin 1970.
\end{thebibliography}
\end{document}